\documentclass{amsart}
\usepackage{latexsym}
\usepackage{amsfonts}
\newtheorem{thm}{Theorem}[section]
 
\newtheorem{lemma}[thm]{Lemma}

\newtheorem{cor}[thm]{Corollary}
\newcommand{\D}{{\mathbf{D}}}
\newcommand{\A}{{\mathbf{A}}}
\newcommand{\PP}{{\mathbf{P}}}
\newcommand{\ds}{\displaystyle\strut} 
\begin{document}
\title{Non-Archimedean Big Picard Theorems} 
\author{William Cherry}
\address{Department of Mathematics, P.O.~Box 311430, University of
North Texas, Denton, TX  76203}
\email{wcherry@unt.edu}
\thanks{Research partially supported by a Junior Faculty 
Research Grant from the University of North Texas}

\date{August 4, 2000}
\commby{}

\begin{abstract} A non-Archimedean analog of the classical Big Picard Theorem,
which says that a holomorphic
map from the punctured disc to a Riemann surface of hyperbolic type extends
accross the puncture, is proven using Berkovich's theory of non-Archimedean
analytic spaces.  
\end{abstract}

\maketitle

\section{Introduction}

One of the fundamental theorems in complex analysis is a theorem of
E.~Picard that says an entire function must take on every complex value
infinitely often, with at most one possible exception.  Today, this is
known as Picard's Little Theorem.  A related theorem of Picard, today
called Picard's Big Theorem, says that in a neigborhood of an isolated
essential singularity, an analytic function must take on every value
infinitely often, again with at most one possible exception.  Picard's
Big Theorem implies his little theorem because a non-polynomial entire
function has an isolated essential singularity at the point at infinity
on the Riemann sphere.

In fact, Picard also proved the following generalization of his little
theorem to Riemann surfaces: Any holomorphic map from the complex plane
to a Riemann surface of genus $\ge 2$ must be constant, and any
holomorphic map from the complex plane to a Riemann surface of genus 1
must hit every point in the curve infinitely often.
An alternative way to phrase Picard's Big Theorem
is to say that a function analytic on a punctured disc which omits at
least two values must extend to a meromorphic function on the whole disc.
Viewed this way, Picard's Big Theorem is an extension theorem, and in
this form it can be generalized to analytic maps between Riemann surfaces.
Recall that a Riemann surface is of hyperbolic type if its holomorphic
universal covering space is biholomorphic to the unit disc.  If a Riemann
surface is not of hyperbolic type, then it is biholomorphic to one
of the following: the Reimann sphere, the complex plane, the complex
plane minus a point, or a complex torus.  Picard proved the following:

\begin{thm}[Picard]\label{picard} Let $Y$ be a Riemann surface
of hyperbolic type,
let $\overline{Y}$ be a compact Riemann surface containing $Y,$
and let $f$ be a holomorphic map from a punctured disc to $Y.$  Then,
$f$ extends to a holomorphic
map from the disc to $\overline{Y}.$
\end{thm}

In \cite{B}, Berkovich proved a non-Archimedean analog of Picard's
Little Theorem, namely any non-Archimedean analytic map from the affine
line to a non-Archimedean Riemann surface of genus $\ge 1$ must be constant.
The purpose of this paper is to show how Berkovich's theory can be used
to prove a non-Archimedean analog of the Big Picard Theorem.

A corollary of Picard's Big Theorem is that if
\hbox{$f:X \rightarrow Y$} is a holomorphic map between 
irreducible algebraic curves
and if $Y$ has genus $\ge 2,$ then $f$ must be a rational map.
This corollary of the Big Picard Theorem was used in a fundamental 
way by Buium in \cite{Bu} to give a new proof of the 
function field version of the following theorem of Raynaud \cite{R}:

\begin{thm}[Raynaud]\label{buim} Let $k\subset F$ be algberaically
closed fields of characteristic zero.  Let $A$ be an Abelian variety
defined over $F$ with $F/k$ trace zero.  Let $X$ be a
projective algebraic curve
in $A$ with geometric genus $\ge 2,$ and let $\Gamma$ be a finite rank
subgroup of $A(F).$  Then, $\Gamma\cap X(F)$ is finite.
\end{thm}

Raynaud did not make use of the Picard theorem.  He used ``reduction
mod $p^2$'' to reduce the theorem to the Mordell conjecture, which was
proven in the characteristic zero function field case by Manin \cite{M}.
Buium's observation that the Picard theorem could be used
to prove Raynaud's Theorem allowed him to generalize the theorem to
higher dimensional $X$ by using higher dimensional
generalizations of Picard's Big Theorem proved by Griffiths and King \cite{GK}
and by Kobayashi and Ochiai \cite{KO}, which was the point of \cite{Bu}.
The theorem Buium proved had
been a conjecture of S.~Lang, see for instance \cite{La},
\cite[Conjecture~I.6.3]{LaBook}.
E.~Hrushovski \cite{H} was able to prove Raynaud's theorem and Buium's 
higher dimensional generalization without the
characteristic zero hypothesis by replacing Buium's use of the Picard
theorem with a model theoretic argument.  To date, this remains the
only proof valid in characteristic $p.$  It would, however, be nice to have
a proof of the characteristic $p$ case of Lang's Conjecture that
does not involve model theory.  One might try to generalize Buium's 
approach to characterstic $p,$ and the first step in doing so would be
to prove a non-Archimedean version of the Big Picard Theorem, which is what
is done in this paper.

\textbf{Acknowledgement.} The author is grateful to A.~Buium for suggesting
that he try to prove a non-Archimedean version of the Big Picard Theorem,
for explaining its potential significance to Lang's Conjecture in
characteristic $p,$ and for his comments on the introduction to this paper.

\section{Basic Notions}

Throughout this paper $k$ will denote fixed algebraically closed field,
complete with respect to a non-trivial non-Archimedean absolute value $|~|.$
The notation $|k|$ will be used to denote the set
\hbox{$\{|z| : z\in k\}.$}  Because $k$ is algebraically closed and
$|~|$ is non-trivial, $|k|$ is dense in the non-negative real numbers.
The quotient
$$
	\tilde k = \{z \in k : |z| \le 1\} / \{z \in k : |z| < 1 \}
$$
is a field known as the \textbf{residue class field} of $k.$

By \textbf{analytic}, I will always mean analytic
in the sense of Berkovich \cite{B}.  Everything in this paper will be
assumed to be defined over $k.$ Furthermore, analytic spaces will always
be assumed to be reduced and separated.  Thus analytic space means
a reduced, separated, Berkovich analytic
space defined over $k,$ and analytic map means a Berkovich analytic map
defined over $k.$

The basic building blocks of non-Archimedean analytic spaces are 
``affinoid spaces.''  The $n$-variable \textbf{Tate Algebra} over $k,$
denoted \hbox{$k{<}z_1, \dots, z_n{>}$}
is the Banach algebra consisting of those formal power series
$\sum_\gamma a_\gamma z^\gamma,$ where the coefficients $a_\gamma$
are in $k$  and
$$
	\lim_{|\gamma|\to\infty} |a_\gamma| = 0.
$$
Here $\gamma$ is a multi-index.  The Banach algebra norm on
\hbox{$k{<}z_1,\dots,z_n{>}$} is the ``$\sup$-norm.'' Namely,
$$
	\left | \sum_\gamma a_\gamma z^\gamma \right |_{\sup}
	= \;\sup_\gamma |a_\gamma|.
$$
The Berkovich analytic space associated to \hbox{$k{<}z_1,\dots,z_n{>}$}
is referred to as the unit $n$-ball and denoted $\mathbf{B}^n$
because the $k$-points of $\mathbf{B}^n$ are those points
$(z_1,\dots,z_n)$ in $k^n$ such that 
\hbox{$\max\{|z_j| : 1\le j \le n\} \le 1.$}  The spaces
$\mathbf{B}^n$ play the role in non-Archimedean analytic geometry that
affine spaces play in algebraic geometry.  If a Banach algebra 
$A$ defined over $k$ is isomorphic to the quotient of a Tate algebra
by a closed ideal, then $A$ is called an \textbf{affinoid algebra.}
The Berkovich analytic space associated to such an algebra is called
an \textbf{affinoid space}.  Note that in Berkovich's terminology
these would be called ``strictly'' affinoid.

Any algebraic variety $X$ over $k$ can be made into a Berkovich analytic
space in a natural way.  I use $X,$ $\PP^1,$ $\A^1,$ \textit{etc{.}} to denote
the analytic space associated to a variety $X,$ the projective line,
the affine line, \textit{etc.}  The $k$-points of these spaces will be denoted
$X(k),$ $\PP^1(k),$ $\A^1(k),$ \textit{etc.}

Given $r>0,$ the ``closed disc'' of radius $r$
and ``center'' $0$ will be denoted $\D(r),$ which is the the affinoid
analytic space associated to the Tate algebra
$$
	k{<}r^{-1}z{>} = \left\{\sum_{n=0}^\infty a_n z^n : 
	\lim_{n\to\infty}|a_n|r^n = 0\right\},
$$
which is also the ring of analytic functions on $\D(r).$  The notation
$\D$ will stand for $\D(1).$  When \hbox{$r\in|k|,$} $\D(r)$ and
$\D$ are clearly isomorphic over $k.$  The notation
$\D^\times(r)$ will be used to denote the ``punctured'' disc of radius
$r,$ which is not affinoid, but whose analytic functions are of the form
$$
	\sum_{n=-\infty}^\infty a_n z^n, \qquad\textrm{with~}
	\lim_{n\to\infty}|a_n|r^n = 0 \textrm{~and~}
	\lim_{n\to-\infty}|a_n|t^n = 0 \textrm{~for all~}t>0.
$$
Occasionally I will need to refer to a disc with non-zero ``center,''
so $\D(a,r)$ will be used to denote the ``closed disc'' of radius
$r$ and ``center'' $a.$  

\begin{thm}[Big Picard]\label{bigpicard} Let $f$ be a meromorphic
function on $\D^\times(r)$ which omits at least two values in $\PP^1(k).$
Then, $f$ extends to an analytic map from $\D(r)$ to 
$\PP^1.$
\end{thm}

\begin{proof}
By shrinking $r$ a little bit if necessary, we may assume
$r$ is in $|k|.$  By post-composition with a M\"obius transformation, it
suffices to assume that $f$ omits the values $0$
and $\infty.$  Thus, we may assume $f$ is an analytic
function on $\D^\times(r)$ and can be written
$$
	f = \sum_{n=-\infty}^\infty a_n z^n,
	\qquad\textrm{with~}\lim_{n\to\infty}|a_n|r^n=0
	\textrm{~and~}\lim_{n\to-\infty}|a_n|t^n = 0
	\textrm{~for all~}t>0.
$$
For $0<t\le r,$ let
$$
	|f|_t = \sup_n \{|a_n| t^n \} = 
	\sup \{|f(z)| : z \in k \textrm{~with~} |z|=t\}.
$$
The last equality is the non-Archimedean maximum modulus principle --
see \cite[Proposition 5.1.3]{BGR}.
The valuation polygon of $f$ is the graph of $\log|f|_t$ as a function
of $\log t.$  The valuation polygon is a piecewise linear graph with
corners at those points $(\log t_0, \log |f|_{t_0})$ where there exist 
two integers $n_1\ne n_2$ such that
$$
	|a_{n_1}|t_0^{n_1} = |f|_{t_0} = |a_{n_2}| t_0^{n_2}.
$$
The theory of Newton or valuation polygons, which essentially
amounts to Hensel's Lemma (or the $p$-adic Weierstrass Perparation Theorem)
and Gauss's Lemma, says that for each corner
$(\log t_0, \log |f|_{t_0})$ of the valuation polygon, $f$ has
a zero $z_0$ with $|z_0|=t_0$ -- see \cite[Ch.~4]{A}.  In fact, the
number of zeros, counting multiplicity, is determined by the 
``sharpness'' of the corner.

If there are infinitely many negative integers $n$ with
$a_n\ne0,$ then it is clear that the valuation polygon for $f$
has infinitely many corners and hence that $f$ has infinitely many
zeros.  Since $f$ is assumed to be zero free, there can be at most
finitely many negative $n$ with $a_n\ne0,$ and hence $f$ has at worst
a pole at $z=0.$ Thus, $f$ extends to a holomorphic map from
$\D(r)$ to $\PP^1.$
\end{proof}

\begin{thm}[Open Mapping]\label{openmapping}
Let \hbox{$r\in |k|$} and let $f$ be a non-constant analytic
function on $\D(r).$  Then 
there exists a $\delta>0$ and in $|k|$ such that the image of $\D(r)$
under $f$ is precisely $\D(f(0),\delta).$
\end{thm}

\begin{proof}
We first consider the $k$-points.  Write
$$
	f(z) - f(0) = \sum_{n=m}^\infty a_n z^n,
$$
with $m\ge1$ and $a_m\ne 0.$  Then for any $w$ in $k,$
$$
	f(z) - w = f(0) - w + \sum_{n=m}^\infty a_n z^n.
$$
From the theory of valuation (Newton) polygons, we see that the series
$$
	f(0) - w + \sum_{n=m}^\infty a_n z^n
$$
will have a zero in $\D(r)(k)$ if and only if
$$
	|f(0) - w | \le \sup_{n\ge m}|a_n|r^n.
$$
Hence, if we let 
$$
	\delta = \sup_{n\ge m}|a_n|r^n, \quad
	\textrm{then}\quad f\big(\D(r)(k)\big) = \D(f(0),\delta)(k).
$$
By \cite[Proposition~2.1.15]{B},
$\D(r)(k)$ is dense in $\D(r)$ and $\D(f(0),\delta)(k)$ is
dense in $\D(f(0),\delta).$
Because $\D(r)$ and $\D(f(0),\delta)$ are compact 
(see \cite[Proposition~1.2.3]{B}) and $f$ is continuous,
the image of $\D(r)$ under $f$ must be $\D(f(0),\delta).$
\end{proof}

The reason that I prefer to work with Berkovich's notion of analytic
spaces is that they are rather nice geometric spaces, unlike the
more traditional ``rigid analytic'' spaces.  For example, we have
the following basic geometric fact.

\begin{thm}\label{connected}The Berkovich analytic space $\D^\times(r)$
is arc-connected and simply connected.
\end{thm}

\begin{proof}
See \cite[\textit{e.g.,} Theorem 4.2.1]{B}.
\end{proof}

The following lemma concludes this section.

\begin{lemma}\label{compose}Let $X$ and $Y$ be connected one dimensional
non-singular analytic spaces, defined over $k.$
Let \hbox{$h\colon X \rightarrow Y$} be a 
non-constant analytic map with finite fibers defined over $k,$ and let 
\hbox{$f\colon\D^\times(r)\rightarrow X$} be an analytic map defined
over $k$ such that
$h\circ f$ extends to an analytic map from $\D(r)$ to $Y.$  Then,
$f$ itself extends to an analytic map from $\D(r)$ to $X.$
\end{lemma}

\begin{proof}
Let $y_0 = h(f(0)),$ which is a point in $Y(k).$
Let $x_j$ be the points in $h^{-1}(y_0),$ which form a finite set
of points in $X(k).$  We can find an affinoid neighborhood $V$ of $y_0$
and affinoid neighborhoods $U_j$ of each point
$x_j$ and each isomorphic to $\D$  such that 
that $U_i\cap U_j$ is empty for all $i\ne j,$ and such that
for each $j,$ $h(U_j)$ is contained in $V.$
Because $h(f(0))=y_0,$ we have 
$h(f(\D(\varepsilon)))$ contained in $V,$ for $\varepsilon$ sufficiently
small.  Because $\D^\times(\varepsilon)$ is connected by 
Theorem~\ref{connected}, this implies
that $f(\D^\times(\varepsilon))$ is contained in just one of the sets $U_j.$ 
Because $U_j$ is isomorphic to $\D,$ we can apply 
Theorem~\ref{bigpicard} to conclude that $f$ extends to an analytic
map from $\D(\varepsilon)$ to $U_j,$ and hence $f$ extends to an analytic
map from $\D(r)$ to $X.$
\end{proof}

\section{Reduction and Uniformization of Curves}

Let $\Gamma$ be a subgroup of $PGL(2,k).$  A point $z$ in $\PP^1$
is called a \textbf{limit point} of $\Gamma$ if there exist
a point $w$ in $\PP^1$ and an infinite sequence $\gamma_n$ in $\Gamma$
such that $\lim \gamma_n(w)=z.$  We let $\Sigma_\Gamma$ denote
the set of limit points of $\Gamma.$  The group $\Gamma$ is called
\textbf{discontinuous} if $\Sigma_\Gamma\ne\PP^1$ and if for each
$z$ in $\PP^1,$ the orbit of $z$ under $\Gamma$ is compact.
The group $\Gamma$ is called 
a \textbf{Schottky group} if it is finitely generated, discontinuous,
and contains no non-trivial torsion elements.

Recall that a set is called \textbf{perfect} if it is equal to its
set of limit points.  In other words a perfect set is a closed set which
does not contain any isolated points.

\begin{thm}\label{mumfordcurve}Let $\Gamma$ be a Schottky group of
rank $g\ge 1$ with
limit set $\Sigma_\Gamma.$  Then, $\Sigma_\Gamma$ is contained
in $\PP^1(k).$  Let
\hbox{$\Omega=\PP^1\setminus\Sigma_\Gamma.$}  Then, $\Gamma$ acts
freely on $\Omega,$ the quotient space $X=\Omega/\Gamma$
is a smooth projective algebraic curve of genus $g,$ and the natural map
\hbox{$\pi:\Omega\rightarrow X$} is an analytic universal covering map.
Moreover, if $g=1,$ then $\Sigma_\Gamma$ consists of exactly two points
and if $g\ge 2,$ then $\Sigma_\Gamma$ is a perfect set.
\end{thm}

\begin{proof}
See \cite{GP} and \cite[\S 4.4]{B}.
\end{proof}

A curve $X$ which is the quotient of a Schottky group as in 
Theorem~\ref{mumfordcurve} is called a \textbf{Mumford curve.}

Given an analytic space, it is possible to associate to it an algebraic
variety $\widetilde{X}$ defined over the residue class field $\tilde k.$
The variety $\widetilde{X}$ is called a \textbf{reduction} of $X$ 
and need not be unique.  I will now recall this notion.

In the case that $X$ is affinoid, there is a canonical reduction
$\widetilde{X}.$ To describe this
more fully, recall that any affinoid algebra $A$ 
has a semi-norm known as the $\sup$-semi-norm and denoted $|~|_{\sup}.$

Thus, one defines
$$
	\widetilde{A} = \{f \in A : |f|_{\sup} \le 1\}
	~/~\{f \in A : |f|_{\sup} < 1 \},
$$
which is an algebra over $\tilde k.$  In fact, since $A$
is affinoid and therefore a quotient of a Tate algebra by a closed ideal,
$\widetilde{A}$ is the quotient of a polynomial ring over 
$\tilde k$ by an ideal, and hence
\hbox{$\widetilde{X} = \mathrm{Spec~} \widetilde{A}$}
is an affine algebraic variety of $\tilde k.$
For example, the reduction of the $n$-ball $\mathbf{B}^n$ is the affine
space $\mathbf{A}^n$ over $\tilde k.$
Affinoid analytic spaces also have a canonical reduction map 
\hbox{$\pi\colon X \rightarrow \widetilde{X}.$}
See \cite[\S 1.3]{B} or \cite[\S 6.2]{BGR} for more details 
on the supremum semi-norm, the reduction of affinoid spaces,
and the reduction map $\pi.$

Before discussing the reduction of more general analytic spaces,
I must recall the notion of formal admissible affinoid coverings.
Let $V$ be an affinoid subspace of an affinoid space $U.$
The inclusion of $V$ into $U$ induces a morphism from
$\widetilde{V}$ to $\widetilde{U}.$  If the induced morphism on
the reductions is an open immersion, then $V$ is called a
\textbf{formal} affinoid subdomain of $U.$  
Recall that by an \textbf{admissible affinoid cover} $\mathcal{U}$ of
an analytic space $X,$ one means a cover $\mathcal{U}$ consisting of
affinoid subspaces $U$ such that if $V$ is any affinoid subspace of $X,$ then
$\mathcal{U}|_V = \{ U \cap V \}$ is a \textit{finite} covering of $V.$
An admissible affinoid covering $\mathcal{U}$ is called \textbf{formal}
if \hbox{$U \cap V$} is a formal affinoid subdomain of $U$ for every
pair $U,V$ in $\mathcal{U}.$  Given a formal admissible affinoid covering 
$\mathcal{U}$ of $X,$ one gets an algebraic variety 
$\widetilde{X}_{\mathcal{U}}$ defined over $\tilde k$ and a
reduction map \hbox{$\pi\colon X \rightarrow \widetilde{X}_\mathcal{U}.$}
In general, different coverings may
give rise to non-isomorphic reductions.

In general, all reductions $\widetilde{X}_\mathcal{U}$ may be singular
varieties even if $X$ is a nice non-singular space.  If 
$\widetilde{X}_\mathcal{U}$ is non-singular, then $X$ is said to have
\textbf{good reduction}.  Fortunately in the
case of curves, one has some control over how bad the singularities
of $\widetilde{X}_\mathcal{U}$ can be.

\begin{thm}[Semi-Stable Reduction]\label{semistable}
If $X$ is a non-singular connected projective
algebraic curve, then there exists a formal admissible affinoid
covering $\mathcal{U}$ of $X$ so that $\widetilde{X}_\mathcal{U}$
is reduced with at worst ordinary double point singularites, and
\begin{enumerate}
\item[(a)] If $X$ has genus $0$ then $\widetilde{X}_\mathcal{U}$ is
isomorphic to $\PP^1.$
\item[(b)] If $X$ has genus $1,$ then 
every non-singular rational component of $X_\mathcal{U}$ meets the other
components in at least two points.  (Such a reduction is called
``semi-stable.'')  Moreover, if $X_\mathcal{U}$ is non-singular,
then up to isomorphism, there exists at most one semi-stable reduction of
$X.$
\item[(c)] If $X$ has genus $\ge 2,$ then every non-singular rational
component of $\widetilde{X}_\mathcal{U}$ meets the other components in
at least three points.  (Such a reduction is called ``stable.'')
Moreover, up to isomorphism, there exists at most one stable reduction
of $X.$
\end{enumerate}
Moreover, $\ds \dim H^1(\widetilde{X}_\mathcal{U},
\mathcal{O}_{\widetilde{X}_\mathcal{U}}) = \dim H^1(X,\mathcal{O}_X) = g.$
\end{thm}

\begin{proof}
See \cite{BL}.
\end{proof}

Note that if $X$ has genus $\ge1$ and $\widetilde{X}_\mathcal{U}$
is a reduction of $X$ assured by Theorem~\ref{semistable},
then either every component of $\widetilde{X}_\mathcal{U}$
is rational or $\widetilde{X}_\mathcal{U}$ has at least one
non-rational component.  This dichotomy  does not depend on the
covering $\mathcal{U}$ because if $X$ has genus~1, the only case where
the specified reduction may not be unique, and if $\widetilde{X}_\mathcal{U}$
is singular, then every component must be rational.  A curve of genus
$\ge 1$ is said to have \textbf{totally degenerate reduction}
if every component of a reduction $\widetilde{X}_\mathcal{U}$
as in Theorem~\ref{semistable} is rational.  The reduction of
curves with totally degenerate reduction is not much use in studying
Picard type questions for maps into $X.$  Fortunately it is precisely
curves with totally degenerate reduction that have a nice uniformization
theory, and that uniformization theory can be used to prove Picard type
theorems in the totally degenerate reduction case.

\begin{thm}\label{berunif}Let $X$ be a non-singular projective algebraic
curve of genus $g\ge1.$ Then, $X$ is a Mumford curve if and only if
$X$ has totally degenerate reduction.
\end{thm}

\begin{proof}
See \cite[Theorem~4.4.1]{B}.
\end{proof}

Finally, we recall that analytic maps lift to normalizations.

\begin{thm}\label{normalization}Let $f\colon X \rightarrow Y$
be an analytic map between reduced analytic spaces.  Let $N(Y)$
denote the analytic subvariety of non-normal points in $Y.$
Assume that $f^{-1}(N(Y))$ is nowhere dense in $X.$  
Let \hbox{$\pi_X\colon \widehat{X}\rightarrow X$} and
\hbox{$\pi_Y\colon \widehat{Y}\rightarrow Y$} be the normalizations of
$X$ and $Y.$  Then, there exists a unique analytic map $\widehat{f}$
from $\widehat{X}$ to $\widehat{Y}$ such that 
\hbox{$\pi_Y\circ \widehat{f}=f\circ\pi_X.$}
\end{thm}

\textit{Poof.} The classical complex analytic argument \cite{GR}
goes through.  See \cite[\S 3.3]{B} and \cite{L} for some details.
$\quad\Box$

\section{Non-Archimedean Analytic Maps to Algebraic Curves}

\begin{thm}\label{disc}Let $f$ be an analytic map from $\D^\times(r)$ to 
an  irreducible projective algebraic curve $X.$  Suppose one of the
following: (a) $X$ has geometric genus $\ge 2;$ (b) $X$ has geometric
genus $1$ and either $f$ omits at least one point of $X$ or
the normalization of $X$ has good reduction; or (c) $X$ has
geometric genus $0$ and $f$ omits at least two points of $X.$
Then $f$ extends to an analytic map from $\D(r)$ to $X.$
\end{thm}

\begin{proof}
By Theorem~\ref{normalization}, it suffices to 
prove the theorem when $X$ is smooth, and the case when $X$
has genus $0$ is already covered by Theorem~\ref{bigpicard}.

When $X$ has positive  genus, the proof breaks up into two cases.

\textit{Case~1.}  Suppose $X$ is a Mumford curve.  Then, by 
Theorem~\ref{mumfordcurve}, there exists a closed set
$\Sigma$ in $\PP^1(k)$ and a universal covering map $\pi$ from
\hbox{$\Omega=\PP^1\setminus\Sigma$} to $X.$  By Theorem~\ref{connected},
$f$ lifts to a map $\hat f$ from $\D^\times(r)$ to $\Omega.$

If $X$ has genus $\ge 2,$ then $\Sigma$ is a perfect subset of
$\PP^1(k)$ by Theorem~\ref{mumfordcurve}.
Since $\Sigma$ contains more than two points, Theorem~\ref{bigpicard}
tells us that $\hat f$ extends to an analytic map from $\D(r)$
to $\PP^1.$  By Theorem~\ref{openmapping}, there exists $\varepsilon>0,$
such that $\D(f(0),\varepsilon)$ is contained in the image of $\hat f.$
Because $\Sigma$ is perfect, this implies $f(0)$ cannot be in $\Sigma,$
and hence $\hat f$ maps to $\Omega.$  Thus, $\pi\circ\hat f$
is an analytic map from $\D(r)$ to $X$ which extends $f.$

If $X$ has genus $1,$ then by assumption $f$ omits a point $x_0$ of $X.$
Moreover, $X$ is a Tate curve, so we may assume $\Sigma$ consists of
the points $0$ and $\infty$ and that $\pi$ is the quotient of the
multiplicative group modulo a rank~1 multiplicative subgroup.
Thus, $\pi^{-1}(x_0)$ accumulates at both $0$ and $\infty.$
Because $\hat f$ omits the two points in $\Sigma,$ it extends to a map from 
$\D(r)$ to $\PP^1$ by Theorem~\ref{bigpicard}.
If $f(0)$ is not zero or infinity, then
$\pi\circ \hat f$ extends $f$ to a map from $\D(r)$ to $X$ and we are done.
But, $f(0)$ cannot be zero or infinity, for in that case
Theorem~\ref{openmapping} would imply that the image of $f$ contained
infinitely many of the points in $\pi^{-1}(x_0),$ a contradiction.

\textit{Case~2.} If $X$ is not a Mumford curve, then by Theorem~\ref{berunif},
$X$ has a formal
affinoid covering $\mathcal{U}$ such that $\widetilde{X_\mathcal{U}}$
contains an irreducible component with positive geometric genus.
As in Berkovich's proof of the little Picard Theorem, 
\cite[Theorem 4.5.1]{B}, there exists a rational function $h$ from 
$X$ to $\PP^1,$ such that $h\circ f$ is bounded in $\A^1.$
Applying Theorem~\ref{bigpicard}, we see that $h\circ f$ extends to an
analytic map from $\D(r)$ to $\PP^1,$ and then applying Lemma~\ref{compose},
we see that $f$ itself extends to an analytic map from $\D(r)$ to
$X.$
\end{proof}

\begin{cor}Let $f\colon X\rightarrow Y$ be an analytic map between
smooth irreducible
algebraic curves defined over $k.$  If $Y$ has geometric genus $\ge 2,$
then $f$ is a rational map.
\end{cor}

\begin{proof}
Let $\overline{X}$ and $\overline{Y}$ be smooth projective
curves containing $X$ and $Y$ as Zariski open subsets.  Then, 
$\overline{X}\setminus X$ is a finite set of isolated points.  The mapping
$f$ extends to an analytic map in a neighborhood of each point in 
$\overline{X}\setminus X,$ by Theorem~\ref{disc}.  Thus, 
$f$ extends to an analytic map from $\overline{X}$ to $\overline{Y}.$
By the GAGA~Theorem, \cite[\S 3.4]{B}, $f$ is rational.
\end{proof}

\vfill\eject
\section{Epilog: July 9, 2002}
I submitted this paper to the 
\textit{Proceedings of the American Mathematical Society} in August~2000.
There was a considerable delay in getting the manuscript refereed, which might
have been partially my fault.  Once finally refereed,
the referee felt that Theorem~\ref{disc},
the original content of this paper, is a ``relatively simple consequence
of the (deep) semistable reduction theorem quoted as Theorem~3.2.''
I do not disagree, but unlike the referee,
I thought a write up of the details of how that
argument goes, as well as my description of a potential application,
might be a useful contribution to the literature.
Because I am currently working on a manuscript which will be much more
general than what is here, I have decided not to try to publish this
paper elsewhere, but will keep the paper available on my web page and
submit it to a preprint archive.

One other comment is in order.  Theorem~\ref{bigpicard}
can be found in:
\begin{quote}
\textsc{M.~van der Put,}
\textit{Essential singularities of rigid analytic functions,}
Nederl.\ Akad.\ Wetensch.\ Indag.\ Math.\ \textbf{43} (1981),423--429.
MR \textbf{83h}:{14021}
\end{quote}
I always suspected I was not the first to make this observation, and shortly
after I submited my paper I found the above reference.  The referee was
also kind enough to point the above reference out to me.  
As also noted by the referee, Theorem~\ref{bigpicard} can now
also be found in the wonderful book:
\begin{quote}
\textsc{A.~Robert,} \textit{A course in $p$-adic analysis,}
Graduate Texts in Mathematics \textbf{198},
Springer-Verlag, 2000. MR \textbf{2001g}:11182
\end{quote}

\vspace{1in}
\end{document}